\documentclass[11pt,english]{article}
\usepackage{amsfonts,amsmath,amsthm,amscd,amssymb,latexsym}

    \usepackage[T2A]{fontenc}
    \usepackage[cp1251]{inputenc}
    \usepackage[ukrainian]{babel}
    \usepackage{graphicx}

\begin{document}

\title{ОДНОСТАЙНА НЕПЕРЕРВНІСТЬ ПО ПРОСТИХ КІНЦЯХ ВІДОБРАЖЕНЬ З УМОВОЮ
НОРМУВАННЯ}

\author{Євген О.\,Севостьянов, Наталія С. Ількевич}

\theoremstyle{plain}
\newtheorem{theorem}{Теорема}[section]
\newtheorem{lemma}{Лема}[section]
\newtheorem{proposition}{Твердження}[section]
\newtheorem{corollary}{Наслідок}[section]
\theoremstyle{definition}

\newtheorem{example}{Приклад}[section]
\newtheorem{remark}{Зауваження}[section]
\newcommand{\keywords}{\textbf{Key words.  }\medskip}
\newcommand{\subjclass}{\textbf{MSC 2000. }\medskip}
\renewcommand{\abstract}{\textbf{Анотація.  }\medskip}
\numberwithin{equation}{section}

\setcounter{section}{0}
\renewcommand{\thesection}{\arabic{section}}
\newcounter{unDef}[section]
\def\theunDef{\thesection.\arabic{unDef}}
\newenvironment{definition}{\refstepcounter{unDef}\trivlist
\item[\hskip \labelsep{\bf Определение \theunDef.}]}%
{\endtrivlist}

\renewcommand{\figurename}{Мал.}

\maketitle

\begin{abstract}
Вивчаються відображення з розгалуженням, які задовольняють деяку
умову спотворення модуля сімей кривих. У ситуації, коли область
визначення відображень є локально зв'язною на своїй межі,
відображена область є регулярною, а мажоранта, яка відповідає за
спотворення модуля сімей кривих, є інтегровною, доведено, що сім'ї
вказаних відображень з одною умовою нормування є одностайно
неперервними в замиканні вихідної області.
\end{abstract}

\medskip {\bf Equicontinuity by the prime ends of mappings with the
normalization condition.} We study branching mappings that satisfy
some condition of distortion of the modulus of families of paths. In
a situation where the definition domain of mappings is locally
connected on its boundary, the mapped domain is regular, and the
majorant responsible for distortion of the modulus of families of
paths is integrable, it is proved that the families of all specified
mappings with one normalization condition are equicontinuous in the
closure of the given domain.

\section{Вступ}
Дану роботу присвячено відображенням з обмеженим і скінченним
спотворенням, які активно вивчаються останнім часом (див., напр.,
\cite{GRY}--\cite{Vu$_1$}). Нещодавно в наших спільних роботах ми
дослідили ситуації, в яких відображення з так званою оберненою
нерівністю Полецького мають неперервне межове продовження, а їх
сім'ї є одностайно неперервними як у внутрішніх, так і межових
точках області, див., напр., \cite{SevSkv$_1$}--\cite{SevSkv$_3$} і
\cite{SSD}. В даній статті ми розглянемо ще один важливий випадок,
коли відображення можуть мати розгалуження, а області мають складну
структуру, при цьому, відображення фіксують принаймні одну точку
області. Зауважимо, що клас відображень з оберненою нерівністю
Полецького включає до себе відображення з обмеженим спотворенням і
відображення зі скінченним спотворенням довжини (див., напр.,
\cite[теорема~3.2]{MRV$_1$}, \cite[теорема~6.7.II]{Ri} і
\cite[теорема~8.5]{MRSY}).

\medskip
Звернемося до означень. Нехай $y_0\in {\Bbb R}^n,$
$0<r_1<r_2<\infty$ і
\begin{equation}\label{eq1**}
A(y_0, r_1,r_2)=\left\{ y\,\in\,{\Bbb R}^n:
r_1<|y-y_0|<r_2\right\}\,.\end{equation}
Для заданих множин $E,$ $F\subset\overline{{\Bbb R}^n}$ і області
$D\subset {\Bbb R}^n$ позначимо через $\Gamma(E,F,D)$ сім'ю всіх
кривих $\gamma:[a,b]\rightarrow \overline{{\Bbb R}^n}$ таких, що
$\gamma(a)\in E,\gamma(b)\in\,F$ і $\gamma(t)\in D$ при $t \in [a,
b].$ Якщо $f:D\rightarrow {\Bbb R}^n$ -- задане відображення,
$y_0\in f(D)$ і $0<r_1<r_2<d_0=\sup\limits_{y\in f(D)}|y-y_0|,$ то
через $\Gamma_f(y_0, r_1, r_2)$ ми позначимо сім'ю всіх кривих
$\gamma$ в області $D$ таких, що $f(\gamma)\in \Gamma(S(y_0, r_1),
S(y_0, r_2), A(y_0,r_1,r_2)).$ Нехай $Q:{\Bbb R}^n\rightarrow [0,
\infty]$ -- вимірна за Лебегом функція.  Будемо говорити, що {\it
$f$ задовольняє обернену нерівність Полецького} в точці $y_0\in
f(D),$ якщо співвідношення
\begin{equation}\label{eq2*A}
M(\Gamma_f(y_0, r_1, r_2))\leqslant \int\limits_{f(D)} Q(y)\cdot
\eta^n (|y-y_0|)\, dm(y)
\end{equation}
виконується для довільної вимірної за Лебегом функції $\eta:
(r_1,r_2)\rightarrow [0,\infty ]$ такій, що
\begin{equation}\label{eqA2}
\int\limits_{r_1}^{r_2}\eta(r)\, dr\geqslant 1\,.
\end{equation}
Зауважимо, що нерівності~(\ref{eq2*A}) добре відомі в теорії
квазірегулярних відображень і виконуються для них при $Q=N(f,
D)\cdot K, $ де $N(f, D)$ -- максимальна кратність відображення в
$D,$ а $K\geqslant 1$ -- деяка стала, яка може бути обчислена як
$K={\rm ess \sup}\, K_O(x, f),$ $K_O(x, f)=\Vert
f^{\,\prime}(x)\Vert^n/|J(x, f)|$ при $J(x, f)\ne 0;$ $K_O(x, f)=1$
при $f^{\,\prime}(x)=0,$ і $K_O(x, f)=\infty$ при
$f^{\,\prime}(x)\ne 0,$ але $J(x, f)=0$ (див., напр.,
\cite[теорема~3.2]{MRV$_1$} або \cite[теорема~6.7.II]{Ri}).
Відображення $f:D\rightarrow {\Bbb R}^n$ називається {\it
дискретним}, якщо прообраз $\{f^{-1}\left(y\right)\}$ кожної точки
$y\,\in\,{\Bbb R}^n$ складається з ізольованих точок, і {\it
відкритим}, якщо образ будь-якої відкритої множини $U\subset D$ є
відкритою множиною в ${\Bbb R}^n.$ Відображення $f$ області $D$ на
$D^{\,\prime}$ називається {\it замкненим}, якщо $f(E)$ є замкненим
в $D^{\,\prime}$ для будь-якої замкненої множини $E\subset D$ (див.,
напр., \cite[розд.~3]{Vu$_1$}). Означення простого кінця, яке
використовується нижче, може бути знайдено в роботі~\cite{IS}, див.
також~\cite{KR$_1$}--\cite{KR$_2$}. Тут і далі $\overline{D}_P$
позначає поповнення області $D$ її простими кінцями, а
$E_D=\overline{D}_P\setminus D$ -- множина всіх простих кінців в
$D.$ Говоримо, що обмежена область $D$ в ${\Bbb R}^n$ {\it
регулярна}, якщо $D$ може бути квазіконформно відображена на область
з локально квазіконформною межею, замикання якої є компактом в
${\Bbb R}^n,$ крім того, кожен простий кінець $P\subset E_D$ є
регулярним. Зауважимо, що у просторі ${\Bbb R}^n$ кожний простий
кінець регулярної області містить ланцюг розрізів з властивістю
$d(\sigma_{m})\rightarrow 0$ при $m\rightarrow\infty,$ і навпаки,
якщо у кінця є вказана властивість, то він -- простий
(див.~\cite[теорема~5.1]{Na}). Крім того, замикання $\overline{D}_P$
регулярної області $D$ є {\it метризовним}, при цьому, якщо
$g:D_0\rightarrow D$ -- квазіконформне відображення області $D_0$ з
локально квазіконформною межею на область $D,$ то для $x, y\in
\overline{D}_P$ покладемо:
\begin{equation}\label{eq1A}
\rho(x, y):=|g^{\,-1}(x)-g^{\,-1}(y)|\,,
\end{equation}
де для $x\in E_D$ елемент $g^{\,-1}(x)$ розуміється як деяка (єдина)
точка межі $D_0,$ коректно визначена з огляду
на~\cite[теорема~4.1]{Na}.

\medskip
Сформулюємо тепер основний результат даної статті. Для цього, для
областей $D, D^{\,\prime}\subset {\Bbb R}^n,$ $n\geqslant 2,$ точок
$a\in D,$ $b\in D^{\,\prime}$ і вимірної за Лебегом функції
$Q:D^{\,\prime}\rightarrow [0, \infty]$ позначимо через ${\frak
S}_{a, b, Q }(D, D^{\,\prime})$ сім'ю всіх відкритих дискретних і
замкнених відображень $f$ області $D$ на $D^{\,\prime},$ що
задовольняють умову~(\ref{eq2*A}) для кожного $y_0\in D^{\,\prime},$
причому $f(a)=b.$ Виконується наступне твердження.

\medskip
\begin{theorem}\label{th2}
{ Припустимо, що область $D$ має слабо плоску межу, жодна із
зв'язних компонент якої не вироджена. Якщо $Q\in L^1(D^{\,\prime}),$
і область $D^{\,\prime}$ є регулярною, то будь-яке $f\in{\frak
S}_{a, b, Q }(D, D^{\,\prime})$ неперервно продовжується до
відображення $\overline{f}:\overline{D}\rightarrow
\overline{D^{\,\prime}}_P,$ причому,
$\overline{f}(\overline{D})=\overline{D^{\,\prime}}_P$ і сім'я
${\frak S}_{a, b, Q }(\overline{D}, \overline{D^{\,\prime}}),$ яка
складається з усіх продовжених відображень
$\overline{f}:\overline{D}\rightarrow \overline{D^{\,\prime}}_P,$
одностайно неперервна в $\overline{D}.$ }
\end{theorem}

\begin{remark}
Зауважимо, що теорему~\ref{th2} можна застосувати для достатньо
широкого спектру областей~$D^{\,\prime}.$ Зокрема, за теоремою
Рімана регулярною областю є будь-яка однозв'язна область, межа якої
містить більше одної точки. Більше того, кожна обмежена скінченно
зв'язна область є конформним образом області, межа якої складається
зі скінченної кількості кіл і ізольованих точок (див., напр.,
\cite[теорема~V.6.2]{Gol}). Оскільки для конформних відображень
ізольовані точки є усувними, то вихідна область може вважатися
регулярною і без вироджених межових компонент.
\end{remark}

\medskip
{\bf 2. Лема про континуум.} Доведення основного результату
ґрунтується на певних властивостях відображень зі збереженням
діаметру прообразу деякого континууму. Наступна лема за деяких інших
припущень на відображення і області, що розглядаються, була доведена
в~\cite[лема~2, пункт 5)]{SevSkv$_1$}, \cite[лема~4.1]{SevSkv$_2$} і
\cite[лема~4.1]{SevSkv$_3$}. Нехай $h$ -- хордальна відстань в
$\overline{{\Bbb R}^n}$ (див., напр., означення~12.1 в~\cite{Va}).

\medskip
\begin{lemma}\label{lem3}
{ Нехай $n\geqslant 2,$ $D$ і $D^{\,\prime}$ -- області в ${\Bbb
R}^n,$ причому, $D$ має слабо плоску межу, жодна компонента
зв'язності якої не вироджується в точку, а область $D^{\,\prime}$ є
регулярною. Нехай також $A$ -- невироджений континуум в
$D^{\,\prime}$ і $\delta>0.$ Припустимо, $f_m$ -- послідовність
відкритих, дискретних і замкнених відображень області $D$ на
$D^{\,\prime}$ з наступною властивістю: для кожного $m=1,2,\ldots$
знайдеться континуум $A_m\subset D,$ $m=1,2,\ldots ,$ такий, що
$f_m(A_m)=A$ і $h(A_m)\geqslant \delta>0.$ Якщо кожне $f_m$
задовольняє співвідношення~(\ref{eq2*A}) для кожного $y_0\in
D^{\,\prime},$ причому, $Q\in L^1(D^{\,\prime}),$ то знайдеться
$\delta_1>0$ таке, що
$$h(A_m,
\partial D)>\delta_1>0\quad \forall\,\, m\in {\Bbb
N}\,.$$
}
\end{lemma}

\begin{proof}
Через компактність простору~$\overline{{\Bbb R}^n}$ межа області $D$
не порожня і є компактом, так що відстань $h(A_m,
\partial D)$ коректно визначена.

\medskip
Проведемо доведення від супротивного. Припустимо, що висновок леми
не є вірним. Тоді для кожного $k\in {\Bbb N}$ знайдеться номер
$m=m_k$ такий, що $h(A_{m_k},
\partial D)<1/k.$ Можна вважати, що послідовність $m_k$ зростає по $k.$
Оскільки $A_{m_k}$ -- компакт, то знайдуться $x_k\in A_{m_k}$ і
$y_k\in
\partial D$ такі, що $h(A_{m_k},
\partial D)=h(x_k, y_k)<1/k$ (див. малюнок~\ref{fig3}).
\begin{figure}[h]
\centerline{\includegraphics[scale=0.45]{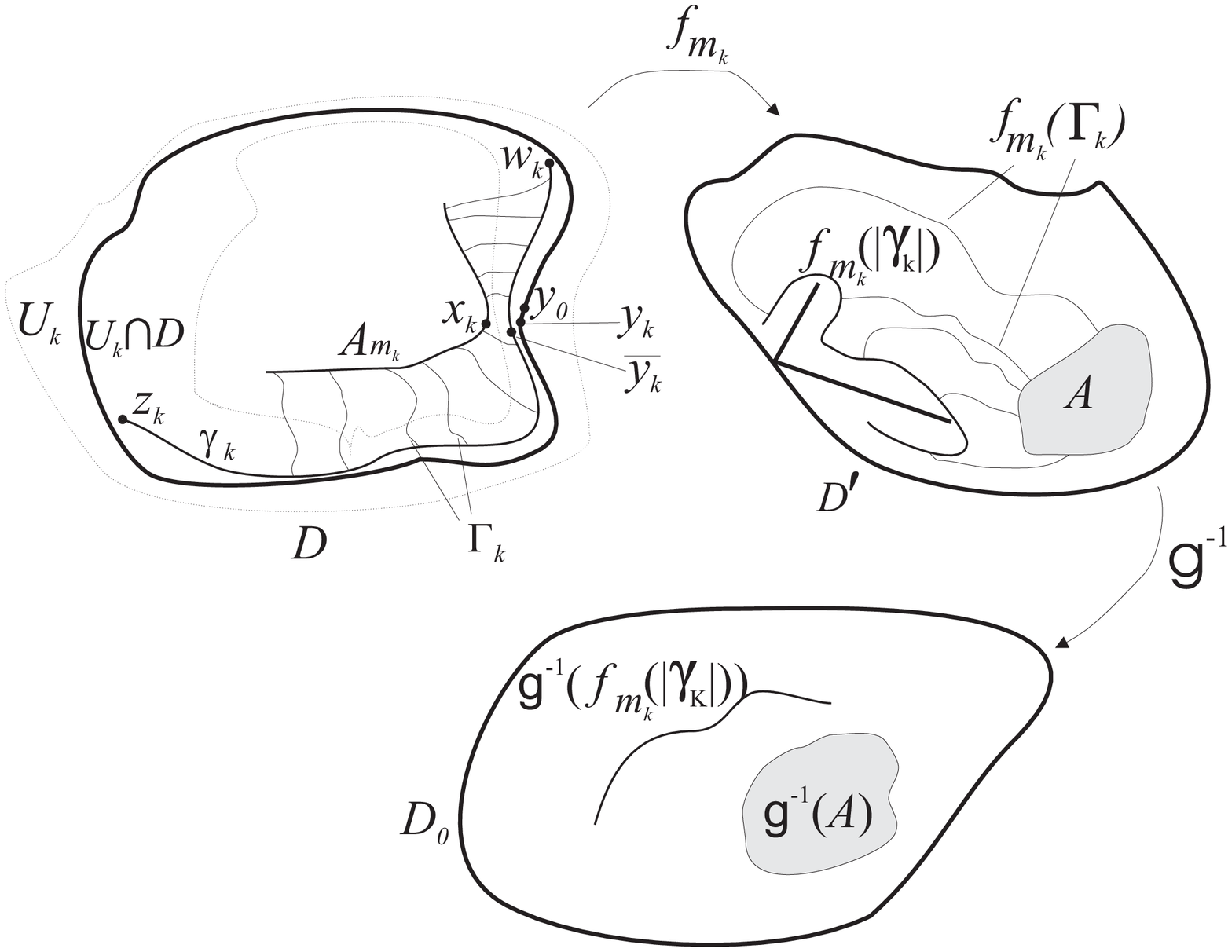}} \caption{До
доведення леми~\ref{lem3}}\label{fig3}
\end{figure}
Оскільки $\partial D$ -- компактна множина, ми можемо вважати, що
$y_k\rightarrow y_0\in
\partial D$ при $k\rightarrow \infty;$ тоді також
%
$x_k\rightarrow y_0\in \partial D$ при $k\rightarrow \infty.$
%
Нехай $K_0$ -- компонента зв'язності $\partial D,$ яка містить точку
$y_0.$ Очевидно, $K_0$ -- континуум в $\overline{{\Bbb R}^n}.$
Оскільки $\partial D$ -- слабо плоска, за теоремою~1 в~\cite{Sev}
відображення $f_{m_k}$ має неперервне продовження
$\overline{f}_{m_k}\colon\overline{D}\rightarrow
\overline{D^{\,\prime}}_P.$ Нехай $\rho$ -- одна з метрик
у~(\ref{eq1A}), і нехай $g:D_0\rightarrow D$ -- квазіконформне
відображення деякої області $D_0$ з локально квазіконформною межею
на $D,$ яке відповідає метриці $\rho$ у~(\ref{eq1A}). Оскільки
$\overline{f}_{m_k}$ неперервне на компакті $\overline{D},$
відображення $\overline{f}_{m_k}$ є рівномірно неперервним у
$\overline{D}$ відносно метрики $\rho$ при кожному фіксованому $k.$
Іншими словами, для кожного $\varepsilon>0$ знайдеться
$\delta_k=\delta_k(\varepsilon)<1/k$ таке, що
\begin{equation}\label{eq3B}
\rho(\overline{f}_{m_k}(x),\overline{f}_{m_k}(x_0))<\varepsilon
\end{equation}
$$
\forall\,\, x,x_0\in \overline{D},\quad h(x, x_0)<\delta_k\,, \quad
\delta_k<1/k\,,$$
де, як звично, $h$ -- хордальна метрика в $\overline{{\Bbb R}^n}.$
Оберемо $\varepsilon>0$ таким, щоб
\begin{equation}\label{eq5D}
\varepsilon<(1/2)\cdot {\rm dist}\,(\partial D_0, g^{\,-1}(A))\,.
\end{equation}
Позначимо $B_h(x_0, r)=\{x\in \overline{{\Bbb R}^n}: h(x, x_0)<r\}.$
Для фіксованого $k\in {\Bbb N},$ покладемо
$$B_k:=\bigcup\limits_{x_0\in K_0}B_h(x_0, \delta_k)\,,\quad k\in {\Bbb
N}\,.$$
Оскільки $B_k$ -- окіл континуума $K_0,$ за~\cite[лема~2.2]{HK}
знайдеться окіл $U_k$ множини $K_0,$ такий, що $U_k\subset B_k$ і
$U_k\cap D$ зв'язна. Можна вважати, що $U_k$ -- відкрита, так що
$U_k\cap D$ є лінійно зв'язною (див.~\cite[пропозиція~13.1]{MRSY}).
Нехай $h(K_0)=m_0.$ Тоді знайдуться $z_0, w_0\in K_0$ такі, що
$h(K_0)=h(z_0, w_0)=m_0.$ Отже, знайдуться послідовності
$\overline{y_k}\in U_k\cap D,$ $z_k\in U_k\cap D$ і $w_k\in U_k\cap
D$ такі, що $z_k\rightarrow z_0,$ $\overline{y_k}\rightarrow y_0$ і
$w_k\rightarrow w_0$ при $k\rightarrow\infty.$ Можна вважати, що
\begin{equation}\label{eq2B}
h(z_k, w_k)>m_0/2\quad \forall\,\, k\in {\Bbb N}\,.
\end{equation}
Оскільки множина $U_k\cap D$ лінійно зв'язна, ми можемо з'єднати
точки $z_k,$ $\overline{y_k}$ і $w_k,$ використовуючи деяку криву
$\gamma_k\in U_k\cap D.$ Як завжди, ми позначаємо через $|\gamma_k|$
носій (образ) кривої $\gamma_k$ в області $D.$ Тоді
$f_{m_k}(|\gamma_k|)$ -- компактна множина в $D^{\,\prime}.$ Якщо
$x\in|\gamma_k|,$ то знайдеться $x_0\in K_0$ таке, що $x\in B(x_0,
\delta_k).$ Зафіксуємо довільне $\omega\in A\subset D.$ Оскільки
$x\in|\gamma_k|$ і, більше того, $x$ -- внутрішня точка $D,$ ми
можемо використовувати запис $f_{m_k}(x)$ замість
$\overline{f}_{m_k}(x).$ Зі співвідношень~(\ref{eq3B})
і~(\ref{eq5D}), а також за нерівністю трикутника, ми отримаємо, що
$$\rho(f_{m_k}(x),\omega)\geqslant
\rho(\omega,
\overline{f}_{m_k}(x_0))-\rho(\overline{f}_{m_k}(x_0),f_{m_k}(x))\geqslant$$
\begin{equation}\label{eq4C}\geqslant {\rm dist}\,(\partial D_0, g^{\,-1}(A))-(1/2)\cdot {\rm
dist}\,(\partial D_0, g^{\,-1}(A))=
\end{equation}
$$=(1/2)\cdot {\rm dist}\,(\partial
D_0, g^{\,-1}(A))>\varepsilon$$
для достатньо великих $k\in {\Bbb N},$ де
$${\rm dist}(\partial D_0, g^{\,-1}(A)):=\inf\limits_{x\in
\partial D_0, y\in g^{\,-1}(A)} |x-y|\,.$$ Переходячи до $\inf$
в~(\ref{eq4C}) по всіх $x\in |\gamma_k|$ і $\omega\in A,$ ми
отримаємо, що
\begin{equation}\label{eq5B}
\rho(f_{m_k}(|\gamma_k|), A)>\varepsilon,\qquad k=1,2,\ldots \,.
\end{equation}

Тепер покажемо, що знайдеться $\varepsilon_1>0$ таке, що
\begin{equation}\label{eq6B}
{\rm dist}\,(f_{m_k}(|\gamma_k|), A)>\varepsilon_1, \quad\forall\,\,
k=1,2,\ldots \,,
\end{equation}
де ${\rm dist},$ як звично, позначає евклідову відстань між
множинами $A, B\subset{\Bbb R}^n.$ Дійсно, нехай
співвідношення~(\ref{eq6B}) не має місця, тоді для числа
$\varepsilon_l=1/l,$ $l=1,2,\ldots$ знайдуться елементи $\xi_l\in
|\gamma_{k_l}|$ і $\zeta_l\in A$ такі, що
\begin{equation}\label{eq7A}
|f_{m_{k_l}}(\xi_l)-\zeta_l|<1/l\,,\quad l=1,2,\ldots \,.
\end{equation}
Можна вважати, що послідовність $k_l,$ $l=1,2,\ldots,$ є зростаючою.
Оскільки~$A$ є компактом, ми також можемо вважати, що послідовність
$\zeta_l$ збігається до $\zeta_0\in A$ при $l\rightarrow\infty.$ За
нерівністю трикутника і по~(\ref{eq7A}) будемо мати, що
\begin{equation}\label{eq8A}
|f_{m_{k_l}}(\xi_l)-\zeta_0|\rightarrow 0\,,\quad
l\rightarrow\infty\,.
\end{equation}
З іншого боку, нагадаємо, що $\rho(f_{m_k}(x),
\omega)=|g^{\,-1}(f_{m_k}(x))-g^{\,-1}(\omega)|,$ де
$g:D_0\rightarrow D$ -- деяке квазіконформне відображення області
$D_0$ на $D,$ див.~(\ref{eq1A}). Зокрема, відображення $g^{\,-1}$ є
неперервним у $D,$ отже, з огляду на нерівність трикутника
і~(\ref{eq8A}), ми отримаємо, що
$$|g^{\,-1}(f_{m_{k_l}}(\xi_l))-g^{\,-1}(\zeta_l)|\leqslant$$
\begin{equation}\label{eq9A}
\leqslant
|g^{\,-1}(f_{m_{k_l}}(\xi_l))-g^{\,-1}(\zeta_0)|+|g^{\,-1}(\zeta_0)-g^{\,-1}(\zeta_l)|\rightarrow
0,\quad l\rightarrow\infty\,.\end{equation}
Проте, за означенням метрики $\rho$ і з~(\ref{eq9A}) випливає, що
$$\rho(f_{m_{k_l}}(|\gamma_{k_l}|), A)\leqslant $$
$$\leqslant \rho(f_{m_{k_l}}(\xi_l), \zeta_l)=
|g^{\,-1}(f_{m_{k_l}}(\xi_l))-g^{\,-1}(\zeta_l)|\rightarrow 0, \quad
l\rightarrow\infty\,,$$
що суперечить~(\ref{eq5B}). Отримана суперечність вказує на
справедливість співвідношення~(\ref{eq6B}).

Покриємо множину $A$ кулями $B(x, \varepsilon_1/4),$ $x\in A.$
Оскільки $A$ компакт, ми можемо вважати, що $A\subset
\bigcup\limits_{i=1}^{M_0}B(x_i, \varepsilon_1/4),$ $x_i\in A,$
$i=1,2,\ldots, M_0,$ $1\leqslant M_0<\infty.$ За означенням, $M_0$
залежить тільки від $A,$ зокрема, $M_0$ не залежить від $k.$
Покладемо
\begin{equation}\label{eq5C}
\Gamma_k:=\Gamma(A_{m_k}, |\gamma_k|, D)\,.
\end{equation}
Нехай $\Gamma_{ki}:=\Gamma_{f_{m_k}}(x_i, \varepsilon_1/4,
\varepsilon_1/2),$ іншими словами, $\Gamma_{ki}$ складається з усіх
кривих $\gamma:[0, 1]\rightarrow D,$ таких що $f_{m_k}(\gamma(0))\in
S(x_i, \varepsilon_1/4),$ $f_{m_k}(\gamma(1))\in S(x_i,
\varepsilon_1/2)$ і $\gamma(t)\in A(x_i, \varepsilon_1/4,
\varepsilon_1/2)$ при $0<t<1.$ Покажемо, що
\begin{equation}\label{eq6C}
\Gamma_k>\bigcup\limits_{i=1}^{M_0}\Gamma_{ki}\,.
\end{equation}
Справді, нехай $\widetilde{\gamma}\in \Gamma_k,$ іншими словами,
$\widetilde{\gamma}:[0, 1]\rightarrow D,$ $\widetilde{\gamma}(0)\in
A_{m_k},$ $\widetilde{\gamma}(1)\in |\gamma_k|$ і
$\widetilde{\gamma}(t)\in D$ при $0\leqslant t\leqslant 1.$ Тоді
$\gamma^{\,*}:=f_{m_k}(\widetilde{\gamma})\in \Gamma(A,
f_{m_k}(|\gamma_k|), D^{\,\prime}).$ Оскільки кулі $B(x_i,
\varepsilon_1/4),$ $1\leqslant i\leqslant M_0,$ утворюють покриття
компакта $A,$ знайдеться $i\in {\Bbb N}$ таке, що
$\gamma^{\,*}(0)\in B(x_i, \varepsilon_1/4)$ і $\gamma^{\,*}(1)\in
f_{m_k}(|\gamma_k|).$ За співвідношенням~(\ref{eq6B}),
$|\gamma^{\,*}|\cap B(x_i, \varepsilon_1/4)\ne\varnothing\ne
|\gamma^{\,*}|\cap (D^{\,\prime}\setminus B(x_i, \varepsilon_1/4)).$
Отже, за~\cite[теорема~1.I.5.46]{Ku} знайдеться $0<t_1<1$ таке, що
$\gamma^{\,*}(t_1)\in S(x_i, \varepsilon_1/4).$ Можна вважати, що
$\gamma^{\,*}(t)\not\in B(x_i, \varepsilon_1/4)$ при $t>t_1.$
Покладемо $\gamma_1:=\gamma^{\,*}|_{[t_1, 1]}.$ З~(\ref{eq6B})
випливає, що $|\gamma_1|\cap B(x_i,
\varepsilon_1/2)\ne\varnothing\ne |\gamma_1|\cap (D\setminus B(x_i,
\varepsilon_1/2)).$ Отже, за~\cite[теорема~1.I.5.46]{Ku} знайдеться
$t_1<t_2<1$ таке, що $\gamma^{\,*}(t_2)\in S(x_i, \varepsilon_1/2).$
Можна вважати, що $\gamma^{\,*}(t)\in B(x_i, \varepsilon_1/2)$ при
всіх $t<t_2.$ Вважаючи $\gamma_2:=\gamma^{\,*}|_{[t_1, t_2]},$
зауважимо, що крива $\gamma_2$ є підкривою $\gamma^{\,*},$ яка
належить $\Gamma(S(x_i, \varepsilon_1/4), S(x_i, \varepsilon_1/2),
A(x_i, \varepsilon_1/4, \varepsilon_1/2)).$

Остаточно, $\widetilde{\gamma}$ має підкриву
$\widetilde{\gamma_2}:=\widetilde{\gamma}|_{[t_1, t_2]},$ таку, що
$f_{m_k}\circ\widetilde{\gamma_2}=\gamma_2,$ причому, $\gamma_2\in
\Gamma(S(x_i, \varepsilon_1/4), S(x_i, \varepsilon_1/2), A(x_i,
\varepsilon_1/4, \varepsilon_1/2)).$ Отже,
співвідношення~(\ref{eq6C}) встановлене. Покладемо
$$\eta(t)= \left\{
\begin{array}{rr}
4/\varepsilon_1, & t\in [\varepsilon_1/4, \varepsilon_1/2],\\
0,  &  t\not\in [\varepsilon_1/4, \varepsilon_1/2]\,.
\end{array}
\right. $$
Зауважимо, що $\eta$ задовольняє співвідношення~(\ref{eqA2}) при
$r_1=\varepsilon_1/4$ і $r_2=\varepsilon_1/2.$ Оскільки відображення
$f_{m_k}$ задовольняє співвідношення~(\ref{eq2*A}), то припускаючи
тут $y_0=x_i,$ отримаємо:
\begin{equation}\label{eq8C}
M(\Gamma_{f_{m_k}}(x_i, \varepsilon_1/4, \varepsilon_1/2))\leqslant
(4/\varepsilon_1)^n\cdot\Vert Q\Vert_1<c<\infty\,,
\end{equation}
де $c$ -- деяка додатна стала і $\Vert Q\Vert_1$ -- $L_1$-норма
функції $Q$ в $D^{\,\prime}.$ З~(\ref{eq6C}) і (\ref{eq8C}),
враховуючи напівадитивність модуля сімей кривих, отримаємо:
\begin{equation}\label{eq4B}
M(\Gamma_k)\leqslant
\frac{4^nM_0}{\varepsilon_1^n}\int\limits_{D^{\,\prime}}Q(y)\,dm(y)\leqslant
c\cdot M_0<\infty\,.
\end{equation}
З іншого боку, оскільки за умовою область $D$ має слабо плоску межу,
з огляду на умову~(\ref{eq2B}), ми отримаємо, що
$M(\Gamma_k)\rightarrow\infty$ при $k\rightarrow\infty,$ що
суперечить~(\ref{eq4B}). Отримане протиріччя доводить лему.
\end{proof}

\section{Доведення теореми~\ref{th2}}

Доведемо теорему~\ref{th2} від супротивного. Припустимо, що ${\frak
S}_{a, b, Q}(\overline{D}, \overline{D^{\,\prime}})$ не є одностайно
неперервною в деякій точці $x_0\in\partial D.$ Тоді знайдуться точки
$x_m\in D$ і відображення $f_m\in {\frak S}_{a, b, Q}(\overline{D},
\overline{D^{\,\prime}}),$ $m=1,2,\ldots ,$ такі що $x_m\rightarrow
x_0$ при $m\rightarrow\infty$ і, причому, при деякому
$\varepsilon_0>0$
\begin{equation}\label{eq15}
h(f_m(x_m), f_m(x_0))\geqslant\varepsilon_0\,,\quad m=1,2,\ldots\,.
\end{equation}
Оберемо довільним чином точку $y_0\in D^{\,\prime},$ $y_0\ne b,$ і
з'єднаємо її з точкою $b$ деякої кривою в $D^{\,\prime},$ яку ми
позначимо $\alpha.$ Покладемо $A:=|\alpha|.$ Нехай $A_m$ -- повне
підняття кривої $\alpha$ при відображенні $f_m$ з початком в точці
$a$ (воно існує за~\cite[лема~3.7]{Vu$_1$}). Зауважимо, що $h(A_m,
\partial D)>0$ за замкненістю відображення~$f_m$ (бо, зокрема, відкриті дискретні
і замкнені відображення є таким, прообраз компакту при яких є
компактом, див.~\cite[теорема~3.3(4)]{Vu$_1$}). Тепер можливі
наступні випадки: або $h(A_m)\rightarrow 0$ при
$m\rightarrow\infty,$ або $h(A_{m_k})\geqslant\delta_0>0$ при
$k\rightarrow\infty$ для деякої зростаючої послідовності номерів
$m_k$ і деякого $\delta_0>0.$

У першому з цих випадків, очевидно, $h(A_m, \partial D)\geqslant
\delta>0$ при деякому $\delta>0.$ Тоді сім'я відображень
$\{f_m\}_{m=1}^{\infty}$ одностайно неперервна в точці $x_0$ за
теоремою~1 в~\cite{Sev}, що суперечить умові~(\ref{eq15}).

У другому випадку, якщо $h(A_{m_k})\geqslant\delta_0>0$ при
$k\rightarrow\infty,$ ми також маємо, що $h(A_{m_k}, \partial
D)\geqslant \delta_1>0$ при деякому $\delta_1>0$ за
лемою~\ref{lem3}. Знову ж таки, теоремою~1 в~\cite{Sev} сім'я
$\{f_{m_k}\}_{k=1}^{\infty}$ є одностайно неперервною в точці $x_0,$
і це суперечить умові~(\ref{eq15}).

Отже, в обох з двох можливих випадків ми прийшли до протиріччя
з~(\ref{eq15}), і це вказує на невірність припущення про відсутність
одностайної неперервності сім'ї ${\frak S}_{a, b, Q}(D,
D^{\,\prime})$ в $\overline{D}.$ Теорема доведена.~$\Box$

\medskip
\begin{example}\label{ex1}
Розглянемо сім'ю відображень $f_n(z)=z^n,$ $n=1,2,\ldots \,,$ $z\in
{\Bbb B}^2=\{z\in {\Bbb C}: |z|<1\}.$ Зауважимо, що $f_n$ є
відображеннями з обмеженим спотворенням як гладкі відображення,
дилатація котрих дорівнює одиниці. Отже, $f_n$ задовольняє
нерівність~(\ref{eq2*A}) при $Q(z)=N(f_n, {\Bbb B}^2)=n,$ де, як
звично, $N$ -- функція кратності, визначена співвідношеннями
$$N(y, f, {\Bbb B}^2)\,=\,{\rm
card}\,\left\{z\in {\Bbb B}^2: f(z)=y\right\}\,, \quad N(f, {\Bbb
B}^2)\,=\,\sup\limits_{y\in{\Bbb C}}\,N(y, f, {\Bbb B}^2)$$
(див., напр., \cite[теорема~3.2]{MRV$_1$} або
\cite[теорема~6.7.II]{Ri}). Всі відображення $f_n$ є дискретними, що
перевіряється безпосередньо, крім того, зберігають межу одиничного
круга і тому є замкненими (див., напр., \cite[теорема~3.3]{Vu$_1$}).
Відображення $f_n$ також фіксують точку $0,$ тому вони задовольняють
всі умови теореми~\ref{th2}, але в той самий час у
нерівності~(\ref{eq2*A}) нема спільної інтегровної функції $Q,$ яка
б забезпечувала всю сім'ю відображень $f_n,$ $n=1,2,\ldots .$
Внаслідок цього, сім'я відображень $f_n$ не є одностайно неперервною
на межі одиничного круга, що перевіряється шляхом безпосередніх
обчислень.
\end{example}

\medskip

\begin{example}
Аналогічний приклад можна побудувати у просторі. Нехай $x\in {\Bbb
B}^n,$ $x=(r\cos\varphi, r\sin\varphi, x_3,x_4,\ldots\,x_n),$ де, як
звично, $x_1=r\cos\varphi,$ $x_2=r\sin\varphi,$ $0\leqslant\varphi<
2\pi,$ $0\leqslant r<\infty.$ Для натурального $m\in {\Bbb N}$
покладемо~$f_m(x)=(r\cos m\varphi, r\sin m\varphi,
x_3,x_4,\ldots\,x_n).$ Шляхом безпосередніх обчислень можна
переконатися, що $K_O(x, f_m)=m^{n-1}$ (див., напр., \cite[приклад~3
пункту~4.3.I]{Re}). Зауважимо, що $f_m$ є відображеннями з обмеженим
спотворенням як гладкі відображення в ${\Bbb R}^n\setminus {\Bbb
R}^{n-2},$ де ${\Bbb R}^{n-2}=\{x\in {\Bbb R}^n: x_{n-1}=x_n=0\},$
дилатація котрих дорівнює $m^{n-1}.$ Отже, $f_m$ задовольняє
нерівність~(\ref{eq2*A}) при $Q(x)=N(f_m, {\Bbb B}^n)\cdot
m^{n-1}=m^n$ в області ${\Bbb B}^n$ (див., напр.,
\cite[теорема~3.2]{MRV$_1$} або \cite[теорема~6.7.II]{Ri}). Всі
відображення $f_m$ є дискретними, що перевіряється безпосередньо,
крім того, зберігають межу одиничної кулі і тому замкнені (див.,
напр., \cite[теорема~3.3]{Vu$_1$}). Відображення $f_m$ також
фіксують точку $0,$ але не мають спільної мажоранти $Q$
в~(\ref{eq2*A}). Неважко переконатися, що сім'я відображень
$\{f_m\}_{m=1}^{\infty}$ не є одностайно неперервною на одиничній
сфері.
\end{example}

КОНТАКТНА ІНФОРМАЦІЯ

\medskip
\noindent{{\bf Євген Олександрович Севостьянов} \\
{\bf 1.} Житомирський державний університет ім.\ І.~Франко\\
кафедра математичного аналізу, вул. Велика Бердичівська, 40 \\
м.~Житомир, Україна, 10 008 \\
{\bf 2.} Інститут прикладної математики і механіки
НАН України, \\
вул.~Добровольського, 1 \\
м.~Слов'янськ, Україна, 84 100\\
e-mail: esevostyanov2009@gmail.com}

\medskip
\noindent{{\bf Наталія Сергіївна Ількевич} \\
Житомирський державний університет ім.\ І.~Франко\\
кафедра фізики і охорони праці, вул. Велика Бердичівська, 40 \\
м.~Житомир, Україна, 10 008 \\ e-mail: ilkevych@1980.gmail.com }


\begin{thebibliography}{99}


\bibitem{GRY} Gutlyanskii~V.~Ya., Ryazanov~V.~I.,
Yakubov~E., \emph{The Beltrami equations and prime ends},
Український математичний вiсник, \textbf{12}, № 1, 2015, p. 27–-66;
translation \emph{The Beltrami equations and prime ends}, J. Math.
Sci. (N.Y.), \textbf{210}, no.~1, 2015, p.~22–-51.

\bibitem{KR$_1$} Ковтонюк Д.А., Рязанов В.И., \emph{К теории простых концов для пространственных
областей}, Укр. мат. журнал, \textbf{67}, № 4, 2015, с. 467--479;
translation ''On the theory of prime ends for space mappings'',
Ukrainian Math. J., \textbf{67}, no.~4, 2015, p.~528–-541.

\bibitem{KR$_2$} Kovtonyuk D.A., Ryazanov V.I., \emph{Prime ends and Orlicz-Sobolev
classes}, St. Petersburg Math. J., \textbf{27}, no.~5, 2016,
p.~765--788.

\bibitem{MRV$_1$} Martio~O., Rickman~S., and V\"{a}is\"{a}l\"{a}~J.,
\emph{Definitions for quasiregular mappings}, Ann. Acad. Sci. Fenn.
Ser. A1., \textbf{448}, 1969, p.~1--40.

\bibitem{MRSY$_1$} Martio O., Ryazanov V., Srebro U. and Yakubov E., \emph{On
$Q$-ho\-me\-o\-mor\-phisms}, Ann. Acad. Sci. Fenn. Ser. A1,
\textbf{30}, no.~1, 2005, p.~49--69.

\bibitem{MRSY} Martio O., Ryazanov V., Srebro U. and Yakubov
E., \emph{Moduli in Modern Mapping Theory.} -- New York: Springer
Science + Business Media, LLC, 2009.

\bibitem{Na} N\"akki R., \emph{Prime ends and quasiconformal
mappings}, J. Anal. Math., \textbf{35}, 1979, p.~13--40.

\bibitem{Re} Reshetnyak~Yu.G., \emph{Space mappings with bounded distortion.} -- Transl.
Math. Monographs 73, AMS, 1989.

\bibitem{Ri} Rickman S., \emph{Quasiregular mappings.} -- Berlin: Springer-Verlag, 1993.

\bibitem{Va} V\"{a}is\"{a}l\"{a} J., \emph{Lectures on $n$-Dimensional
Quasiconformal Mappings}, Lecture Notes in Math. \textbf{229},
Berlin etc.: Springer--Verlag, 1971.

\bibitem{Vu$_1$} Vuorinen~M., \emph{Exceptional sets and boundary behavior of quasiregular
mappings in $n$-space}, Ann. Acad. Sci. Fenn. Ser. A 1. Math.
Dissertationes, \textbf{11}, 1976, p.~1--44.


\bibitem{SevSkv$_1$} Севостьянов~Е.А., Скворцов~С.А.,
\emph{О сходимости отображений в метрических пространствах с прямыми
и обратными модульными условиями}, Укр. мат. журнал, \textbf{70},
№~7, 2018, с.~952--967; transl. Sevost'yanov~E.A., Skvortsov~S.A.,
\emph{On the Convergence of Mappings in Metric Spaces with Direct
and Inverse Modulus Conditions}, Ukr. Math. J., \textbf{70,} no.~7,
2018, p.~1097--1114.

\bibitem{SevSkv$_2$} Севостьянов~Е.А., Скворцов~С.А.,
\emph{О локальном поведении одного класса обратных отображений,}
Укр. мат. вестник, \textbf{15}, №~3, 2018, с.~399--417; translation
Sevost'yanov~E.A., Skvortsov~S.A., \emph{On the local behavior of a
class of inverse mappings}, J. Math. Sci., \textbf{241}, no.~1,
2019, p.~77--89.

\bibitem{SevSkv$_3$} Sevost’yanov, E.A. and Skvortsov~S.A.,
\emph{On mappings whose inverse satisfy the Poletsky inequality},
Ann. Acad. Scie. Fenn. Math., \textbf{45}, 2020, p.~259--277.

\bibitem{SSD} Sevost'yanov~E.A., Skvortsov~S.O., Dovhopiatyi~O.P.,
\emph{On mappings satisfying the inverse Poletsky inequality} //
www. arxiv. org, arXiv:1904.01513.

\bibitem{IS} Ильютко~Д.П., Севостьянов~Е.А., \emph{О простых концах на римановых
многообразиях}, Укр. мат. вестник, \textbf{15}, № 3, 2018,
c.~358--392; translation \emph{On prime ends on Riemannian
manifolds} in J. Math. Sci., \textbf{241}, no.~1, 2019, p.~47--63.

\bibitem{Gol} Голузин Г.М., \emph{Геометрическая теория функций комплексного
переменного.} -- М.: Наука, ФИЗМАТГИЗ, 1966.

\bibitem{Sev} Севостьянов~Е.А., \emph{Межове продовження відображень з оберненою нерівністю
Полецького по простих кінцях}, Укр. мат. журнал (подана до друку).

\bibitem{HK} Herron~J. and Koskela~P., \emph{Quasiextremal distance domains
and conformal mappings onto circle domains}, Compl. Var. Theor.
Appl., \textbf{15}, 1990, p.~167--179.

\bibitem{Ku} Куратовский~К., \emph{Топология}, т.~2. -- М.:
Мир, 1969.

\end{thebibliography}
\end{document}